\documentclass{article}

\usepackage{hyperref}
\usepackage{mathptmx}       
\usepackage{helvet}         
\usepackage{courier}        
\usepackage{type1cm}        
\usepackage{mathtools}
\usepackage{makeidx}         
\usepackage{graphicx}        
\usepackage{multicol}        
\usepackage[bottom]{footmisc}
\usepackage{enumerate}

\makeindex             

\usepackage{hyperref}
\usepackage{color}
\usepackage{mathrsfs}
\usepackage{amsmath}
\usepackage{amssymb}

\newtheorem{theorem}{Theorem}[section]
\newtheorem{lemma}[theorem]{Lemma}

\newtheorem{proposition}[theorem]{Proposition}

\newtheorem{remark}[theorem]{Remark}
\newtheorem{definition}[theorem]{Definition}

\begin{document}

\title{Endpoint results for Fourier integral operators on noncompact symmetric spaces}

\author{Tommaso Bruno\footnote{Dipartimento di Scienze Matematiche ``Giuseppe Luigi Lagrange'', Politecnico di Torino, Corso Duca degli Abruzzi 24, 10129 Torino, Italy, tommaso.bruno@polito.it}, Anita Tabacco \footnote{Dipartimento di Scienze Matematiche ``Giuseppe Luigi Lagrange'', Politecnico di Torino, Corso Duca degli Abruzzi 24, 10129 Torino, Italy, anita.tabacco@polito.it} and Maria Vallarino \footnote{Dipartimento di Scienze Matematiche ``Giuseppe Luigi Lagrange'', Politecnico di Torino, Corso Duca degli Abruzzi 24, 10129 Torino, Italy, maria.vallarino@polito.it}  }

%
%
\maketitle

\abstract{Let $\mathbb X$ be a noncompact symmetric space of rank one and let $\mathfrak h^1(\mathbb X)$ be a local atomic Hardy space. 
We prove the boundedness  from $\mathfrak{h}^1(\mathbb{X})$ to $L^1(\mathbb{X})$ and on $\mathfrak{h}^1(\mathbb{X})$ of some classes of Fourier integral operators related to the wave equation associated with the Laplacian on $\mathbb X$ and we estimate the growth of their norms depending on time. 

\vspace{.5cm}


}

\section{Introduction}

Given a second order differential operator $\mathcal L$ on a manifold $\mathbb M$ consider the Cauchy problem for the associated wave equation 
\begin{equation}\label{PC}
\begin{cases}
&\partial_t^2u(t,x)+\mathcal Lu(t,x)=0,\\
&u(0,x)=f(x),\\
&\partial_tu(0,x)=g(x)\qquad t\in\mathbb R, \;x\in \mathbb M.
\end{cases}
\end{equation}
An interesting problem is to find $L^p$-bounds of the solution $u$ at a certain time in terms of Sobolev norms of the initial data $f$ and $g$. This problem is well understood for the standard Laplacian in $\mathbb R^n$~\cite{Mi,P}. It was also studied for the Laplace--Beltrami operator on compact manifolds~\cite{SSS}, for the subLaplacian on groups of Heisenberg type~\cite{MS, MSt} and for the Laplacian on compact Lie groups~\cite{CFS}. Ionescu~\cite{I} investigated the same problem on noncompact symmetric spaces of rank one. More precisely, let $\mathbb X$ be a noncompact symmetric space of rank one and dimension $n$ and denote by $d$ the number $(n-1)/2$.  Let $\Delta$ denote the Laplace--Beltrami operator on $\mathbb X$, whose $L^2$-spectrum is the half-line $[\rho^2,\infty)$, and set $\mathcal L=\Delta-\rho^2$ (see Section~\ref{notation} for the definition of $\rho$). The wave equation associated with $\mathcal L$ was considered in~\cite{AMPS, APV1, APV2, He2, I, Ta}. By the spectral theorem the solution of the Cauchy problem~\eqref{PC} associated with $\mathcal L$ is given by 
\[
u(t,\cdot)=\cos(t\sqrt \mathcal L)f+\frac{\sin( t\sqrt \mathcal L)}{\sqrt \mathcal L}g.
\]
Finding $L^p$-bounds for $u$ amounts to prove the boundedness on $L^p(\mathbb X)$ of the operators
\[\mathcal T_{t}=m(\sqrt \mathcal L)\cos ({t\sqrt \mathcal L})\qquad {\rm{and}}\qquad \mathcal S_{t}=m(\sqrt \mathcal L)\frac{\sin( t\sqrt \mathcal L)}{\sqrt \mathcal L},\]  
for suitable symbols $m$, and estimate the growth of their norm on $L^p(\mathbb X)$ depending on $t$. 

In this paper we prove endpoint results at $p=1$ for $\mathcal T_{t}$. To state our result, we need some notation. For every $a\geq 0$ and $b\in\mathbb R$ let $S^b_a$ be the set of continuous functions $m$ on the complex tube $\{\lambda\in \mathbb C\colon |\mathrm{Im}\, \lambda|\leq a\}$, analytic in the interior of the tube, infinitely differentiable on the two lines $|\mathrm{Im}\, \lambda|= a$, which satisfy the symbol inequalities
\[
|\partial^{\alpha}_{\lambda}m(\lambda)|\leq C\,(1+|\mathrm{Re}\,\lambda|)^{b-\alpha}\qquad \forall\, \alpha\in\mathbb N,\;\; |\mathrm{Im}\, \lambda|\leq a.
\]
If $m\in S^b_a$, the real number $b$ is called the \emph{order} of $m$.

In~\cite{I}, Ionescu proved an endpoint result for $\mathcal{T}_t$ at $p=\infty$. Indeed, he showed that if $m\in S^{-d}_{\rho}$ is an even symbol, then the operator $\mathcal T_t$ is bounded from $L^\infty(\mathbb{X})$ to a suitable $BMO(\mathbb X)$ space. From this, he deduced the boundedness of $\mathcal{T}_t$ on $L^p(\mathbb{X})$ for every $p\in(1,\infty)$. Let us also mention that previously Giulini and Meda~\cite{GM} proved $L^p$-estimates, $p\in (1,\infty)$, for oscillating multipliers of the form $\Delta^{-\beta/2}\,e^{i\Delta^{\alpha/2}}$, $\alpha>0$, $\mathrm{Re}\,\beta\geq 0$. When $\alpha =1$ and $\beta=d$, these operators are related to $\mathcal T_{1}$. Note, however, that on a noncompact symmetric space the growth in $t$ of the norm of $\mathcal T_{t}$ cannot be deduced from its norm at $t=1$, as one can do in other contexts equipped with a dilation structure (e.g.\ Euclidean spaces and stratified nilpotent groups).

Let $\mathfrak h^1(\mathbb X)$ be the local atomic Hardy space of Goldberg type defined by Taylor~\cite{T} and Meda and Volpi~\cite{MV} (see Definition~\ref{def:HARDY} below). 
The main result of this paper is the following. 

\begin{theorem}\label{main_teo_intro}
Let $t>0$. Then the following hold:
\begin{enumerate}
\item[(i)] if $m\in S^{-d}_{\rho}$ is an even symbol, then the operator $\mathcal T_{t}$ is bounded from $\mathfrak{h}^1(\mathbb X)$ to $L^1(\mathbb{X})$ and $\|\mathcal T_{t}f\|_{\mathfrak h^1\rightarrow L^1}\leq C\,e^{\rho \,t}$;
\item[(ii)] if $m\in S^{b}_{\rho}$ is an even symbol and $b<-d$, then the operator $\mathcal T_{t}$ is bounded on $\mathfrak{h}^1(\mathbb X)$ and $\|\mathcal T_{t}f\|_{\mathfrak h^1\rightarrow \mathfrak h^1}\leq C\,e^{\rho \,t}$.
\end{enumerate}
\end{theorem} 
The results of Theorem~\ref{main_teo_intro} are endpoint results for $\mathcal T_{t}$ at $p=1$. The $\mathfrak{h}^1 \to L^1 $ boundedness may be considered as the counterpart at $p=1$ of Ionescu's result; observe, however, that it does not descend from this by duality, for $\mathfrak{h}^1(\mathbb{X})$ is not the dual of $BMO(\mathbb{X})$. Nevertheless, the proof of part (i) is strongly related to Ionescu's proof. Part (ii), instead, gives a more precise endpoint result but requires higher regularity of the multiplier $m$. It would be interesting to know whether this regularity condition is really necessary, or whether it can be weakened up to the value $b=-d$, which as part (i) shows is enough for the $\mathfrak{h}^1\to L^1$ boundedness. The proof of part (ii) goes through a pointwise decomposition of the convolution kernel $k_{t}$ of $\mathcal{T}_t$ as a sum of compactly supported functions in certain annuli, whose $\mathfrak{h}^1$-norm we estimate separately. We do this by means of precise estimates of both $k_{t}$ and its derivative. In applying this procedure, the condition $b<-d$ turns out to be fundamental.

We finally observe that, by analytic interpolation with $L^2(\mathbb{X})$ and by duality, one can re-obtain Ionescu's result of $L^p(\mathbb{X})$-boundedness of $\mathcal{T}_t$, $p\in (1,\infty)$.
 
\smallskip

The paper is organized as follows. In Section~\ref{notation} we summarize the notation for noncompact symmetric spaces of rank one and the spherical analysis on them. In Section~\ref{h1} we recall the definition of the local Hardy space $\mathfrak h^1(\mathbb X)$ and we prove some technical lemmata which will be of use later on. In Section~\ref{main1} we prove Theorem~\ref{main_teo_intro}~(i), while Section~\ref{main} is devoted to the proof of Theorem~\ref{main_teo_intro}~(ii).
\section{Notation}\label{notation}

We shall use the same notation as in~\cite{I} and refer the reader to~\cite{AJi, GV, H} for more details on noncompact symmetric spaces and spherical analysis on them. 

Let $G$ be a connected noncompact semisimple Lie group with finite centre, $\mathfrak g$ its Lie algebra, $\theta$ a Cartan involution of $\mathfrak g$ and $\mathfrak g=\mathfrak k\oplus \mathfrak p$ the associated Cartan decomposition. Let $K$ be a maximal compact subgroup of $G$ and $\mathbb X=G/K$ be the associated symmetric space of dimension $n$. Let $\mathfrak a$ be a maximal abelian subspace of $\mathfrak p$. We will assume that the dimension of $\mathfrak a$ is one, i.e. that the rank of $\mathbb X$ is one. The Killing form on $\mathfrak g$ induces a $G$-invariant distance on $\mathbb X$, which we shall denote by $d(\cdot,\cdot)$. For every $x\in\mathbb X$ we denote by $|x|$ the distance $d(x,o)$, where $o=e K$ and $e$ is the identity of $G$. Let $\mathfrak a^*$ be the real dual of $\mathfrak a$ and for $\alpha\in\mathfrak a^*$ let $\mathfrak g_{\alpha}=\{X\in\mathfrak g\colon [H,X]=\alpha(H)X \quad \forall\, H\in\mathfrak a\}$. Let $\Sigma=\{\alpha\in\mathfrak a^*\setminus\{0\}\colon {\rm{dim}\,} \mathfrak g_{\alpha}\neq 0\}$ be the set of non-zero roots. It is well known that either $\Sigma=\{-\alpha,\alpha\}$ or $\Sigma=\{-2\alpha,-\alpha,\alpha,2\alpha\}$. Let $m_1={\rm{dim}\,} \mathfrak g_{\alpha}$, $m_2={\rm{dim}\,} \mathfrak g_{2\alpha}$ and $\rho=(m_1+2m_2)/2 $. Set $\mathfrak n=\mathfrak g_{\alpha} +\mathfrak g_{2\alpha}$ and $N={\rm{exp}}\,\mathfrak n$.

In the sequel we shall identify $A={\rm{exp}}\,\mathfrak a$ with $\mathbb R$ by choosing the unique element $H_0$ of $\mathfrak a$ such that $\alpha(H_0)=1$ and considering the diffeomeorphism $a:\mathbb R\rightarrow A$ defined by $a(s)={\rm{exp}}(sH_0)$. It is well known that $G$ admits the Cartan decomposition $G=KA^{+}K$, where $A^+=\{a(s):s\geq 0\}$ and the Iwasawa decomposition $G=NAK$. For every $g\in G$ we denote by $H(g)$ the unique element in $\mathbb R$ such that $g=n\, {\rm{exp}}(H(g)H_0)k$, for some $n\in N$ and $k\in K$. 

For every $r>0$ and $x\in\mathbb X$ we denote by $B(x,r)$ the closed ball centred at the point $x$ of radius $r$. For every $0<r<R$ we denote by $A_r^R$ the annulus $A_{r}^{R}=\{x\in\mathbb X\colon r\leq |x|\leq R\}$. As a convention, $A^R_r$ when $r\leq 0$ shall be intended as the ball $B(o,R)$.

For every integrable function $f$ on $G$ we have
\[
\int_Gf(g)\, dg=C\int_{K}\int_{\mathbb R^+}\int_K f(k_1a(s)k_2) \,\delta(s)\, dk_1 \,ds \,dk_2,
\] 
where $dg$ is the Haar measure of $G$, $dk$ is the Haar measure of $K$ normalized in such a way that $\int_K	\, dk=1$ and 
\[
\delta(s)=C(\sinh s)^{m_1}(\sinh 2s)^{m_2} {\asymp} \,\begin{cases}
s^{n-1}&s\leq 1\\
e^{2\rho s}&s>1.
\end{cases}
\]
We identify right $K$-invariant functions on $G$ with functions on $\mathbb X$, and $K$-biinvariant functions on $G$ with $K$-invariant functions on $\mathbb X$ which can also be identified with functions depending only on the coordinate $s\in \mathbb R^+$. More precisely, if $f$ is a $K$-biinvariant function on $G$ we shall denote by $F:\mathbb R^+\rightarrow \mathbb C$ the function such that $f(k_1a(s)k_2)=F(s)$ for every $s\in\mathbb R^+$, $k_1,k_2\in K$. We define the convolution of two functions $f_1,f_2$ on $\mathbb X$, when it exists, as
\[
f_1\ast f_2(x)=\int_G f_1(gh)f_2(h^{-1})\, dh\qquad \forall\, x=g K\in\mathbb X.
\]
We denote by $\mu$ the Riemannian measure on $\mathbb X$ and for every $p\in [1,\infty)$ let $L^p(\mathbb X)$ be the space of measurable functions $f$ such that $\|f\|_{L^p}^p=\int_{\mathbb X}|f|^pd\mu<\infty$. For every $K$-invariant function $f$ on $\mathbb X$  
\[
\int_{\mathbb X} f(x)\, d\mu(x)=\int_{\mathbb R^+} F(s) \,\delta(s) \,ds,
\]
where $F$ is defined above. By this and the left-invariance of the metric
\begin{equation}\label{balls}
\mu(B(x,r))=\mu(B(o,r)) \asymp \begin{cases}
r^n&r\leq 1\\
e^{2\rho r}&r>1
\end{cases}\qquad \forall \, r>0,\; x\in\mathbb X.
\end{equation}
Observe moreover that 
\begin{equation}\label{misura_anello}
\mu(A^{R+r}_{R-r}) \asymp e^{2\rho R}r, \qquad\forall \, R>1, \; r<1.
\end{equation}
We recall that a spherical Fourier transform on the symmetric space is defined. It associates to each left $K$-invariant function $f$ on $\mathbb X$, i.e. to each radial function, its spherical Fourier transform $\widetilde f$, defined by
\[
\widetilde f(\lambda)=\int_G f(g)\,\phi_{\lambda}(g) \,dg\qquad \lambda\in\mathfrak a^*_{\mathbb C},
\]
where the spherical functions are defined by
\[
\phi_{\lambda}(g)=\int_K{\rm{exp}}[(i\lambda+\rho)H(kg)]\,dk\qquad g\in G,\;\lambda\in\mathfrak a^*_{\mathbb C}.
\] 
It is well known that for every radial function in $L^2(\mathbb X)$ 
\begin{equation}\label{plancherel}
\|f\|_{L^2}^2=C\,\int_0^{\infty}|\widetilde f(\lambda)|^2\,|{\bf{c}}(\lambda)|^{-2}d\lambda,
\end{equation}
and
\begin{equation}\label{inversion}
f(x)=C\,\int_0^{\infty}\widetilde{f}(\lambda)\,\phi_{\lambda}(x)\,|{\bf{c}}(\lambda)|^{-2}d\lambda,
\end{equation}
where ${\bf{c}}$ is the Harish-Chandra function. In particular, by the Plancherel and the inversion formulae above, any bounded function $m:\mathbb R^+\rightarrow \mathbb C$ defines a bounded operator on $L^2(\mathbb X)$ given by $\widetilde{\mathcal U_mf}(\lambda)=m(\lambda)\,\widetilde f(\lambda)$. 
    
All throughout  the paper, we shall write $A\lesssim B$ when there exists a positive constant $C$ such that $A\leq C \, B$, whose value may change from line to line. If $A\lesssim B$ and $B\lesssim A$, we write $A\asymp B$.

\section{The local Hardy space $\mathfrak{h}^1(\mathbb X)$}\label{h1}
We recall here the definition of the local atomic Hardy space $\mathfrak h^1(\mathbb X)$, which can be thought as the analog in the context of noncompact symmetric space of the local Hardy space introduced by Goldberg in the Euclidean setting~\cite{G}. The space $\mathfrak h^1(\mathbb X)$ was introduced and studied by Meda and Volpi~\cite{MV} and Taylor~\cite{T} in more general contexts. It is easy to see that noncompact symmetric spaces satisfy the geometric assumptions of~\cite{MV} and~\cite{T}, so that the theory developed in those papers can be applied in our setting.  

\begin{definition}
A {\bf{standard $\mathfrak{h}^1$-atom}} is a function $a$ in $L^1(\mathbb X)$ supported in a ball $B$ of radius $\leq 1$ such that
\begin{enumerate}[(i)]
\item $\|a\|_{L^2}\leq \mu(B)^{-1/2}$ (size condition);
\item $\int a \,d\mu=0$ (cancellation condition).   
\end{enumerate}
A {\bf{global $\mathfrak{h}^1$-atom}} is a function $a$ in $L^1(\mathbb X)$ supported in a ball $B$ of radius $1$ such that $\|a\|_2\leq \mu(B)^{-1/2}$. Standard and global $\mathfrak{h}^1$-atoms will be referred to as {\bf{admissible atoms}}.  
\end{definition}
\begin{definition}\label{def:HARDY}
The Hardy space $\mathfrak{h}^1(\mathbb X)$ is the space of functions $f$ in $L^1(\mathbb X)$ such that $f=\sum_jc_ja_j$, where $\sum_j|c_j|<\infty$ and $a_j$ are admissible atoms. The norm $\|f\|_{\mathfrak{h}^1}$ is defined as the infimum of $\sum_j|c_j|<\infty$ over all atomic decompositions of $f$. 
\end{definition}

By means of the atomic structure of $\mathfrak{h}^1(\mathbb{X})$ and of the following result, the boundedness from $\mathfrak{h}^1(\mathbb{X})$ of an operator bounded on $L^2(\mathbb X)$ may be tested only on atoms.  Its proof is an easy adaptation of the proof of~\cite[Theorem 4 and Proposition 4]{MV} and is omitted.
\begin{proposition}\label{uniform}
Let $Y$ be either $L^1(\mathbb{X})$ or $\mathfrak{h}^1(\mathbb{X})$. Suppose that $\mathcal U$ is a $Y$-valued linear operator defined on finite linear combination of admissible atoms such that
\[
A\coloneqq \sup\{\|\mathcal Ua\|_{Y} \colon a \mbox{ $\mathfrak{h}^1$-atom} \}<\infty.
\] 
Then there exists a unique bounded operator $\mathcal{U}'$ from $\mathfrak h^1(\mathbb X)$ to $Y$ which extends~$\mathcal U$ with norm $\| \mathcal{U}'\|_{\mathfrak{h}^1\to Y} \lesssim A$. If $\mathcal{U}$ is bounded on $L^2(\mathbb X)$, then $\mathcal{U}'$ and $\mathcal U$ coincide on $Y\cap L^2(\mathbb X)$.
\end{proposition}

\bigskip

We now collect some technical lemmata where we estimate the $\mathfrak h^1$-norm of $L^2$-functions supported either in a ball or in an annulus, which will be useful later on. We shall repeatedly use the notion of discretization of the space $\mathbb X$, which we now recall.

For every $r\in (0,1]$, we call $r/3$-discretization $\Sigma$ of $\mathbb X$ a set of points which is maximal with respect to the properties
\[
\min\{d(z,w):z,w\in\Sigma, z\neq w\}>\frac{r}{3}, \quad \qquad d(x,\Sigma)\leq \frac{r}{3}\quad \forall\,  x\in \mathbb X.
\]
Let  $\Sigma$ be a $r/3$-discretization of $\mathbb X$, for some $r\in (0,1]$. Then the family of balls $\mathcal B=\{B(z,r):z\in \Sigma\}$ is a uniformly locally finite covering of $\mathbb X$. More precisely, there exists a constant $M$, independent of $r$, such that
\begin{equation}\label{Muniform}
1\leq \sum_{B\in\mathcal B}\chi_B(x)\leq M\qquad \forall \, x\in\mathbb X.
\end{equation}
Indeed, given any point $x\in \mathbb X$, if $x\in B(z,r)$, then $z\in B(x,r)$. Thus $\sum_{B\in\mathcal B}\chi_B(x)=M(x)=|\Sigma\cap B(x,r)|$. Let $\{w_1,\dots, w_{M(x)}\}=\Sigma\cap B(x,r)$. If $w_i,w_j\in \Sigma\cap B(x,r)$, with $w_i\neq w_j$, then 
$B(w_i,\frac{r}{6})\cap B(w_j,\frac{r}{6})=\emptyset$. Thus $\bigcup_{i=1}^{M(x)} B(w_i,\frac{r}{6})\subseteq B(x,r+\frac{r}{6})$ and by~\eqref{balls}
\[
C\,M(x) r^n\leq \mu\left(  \bigcup_{i=1}^{M(x)} B(w_i,\tfrac{r}{6})  \right) \leq \mu\left( B(x,r+\tfrac{r}{6}) \right)\leq C \,r^n.
\]
Thus there exists a constant $M$ independent of $x$ and $r$ such that $M(x)\leq M,$ which proves~\eqref{Muniform}.
 
 \begin{lemma}\label{normah1funzioneL2}
Let $f$ be a function in $L^2(\mathbb X)$ supported in a ball $B=B(o,R)$. If
\begin{itemize}
\item either $ R\leq 1$ and $f$ has vanishing integral,
\item or $R\geq 1$,
\end{itemize}
then $\|f\|_{\mathfrak{h}^1}\lesssim  \mu(B)^{1/2}\,\|f\|_{L^2}$.
\end{lemma}
{\bf{Proof.}} If $ R\leq 1$ and $f$ has vanishing integral, it suffices to notice that $\frac{f}{\mu(B)^{1/2}\|f\|_{L^2}}$ is a standard atom. 

If $R\geq 1$, we follow the line of~\cite[Lemma 3.3]{MV} with slight modifications. Let $\Sigma$ be a $1/3$-discretization of $\mathbb X$. Denote by $z_1,\dots,z_N$ the points in $ \Sigma$ such that $B(z_j,1)\cap B\neq \emptyset$. Note that $N\leq C\, \mu(B)$. Denote by $B_j$ the ball $B(z_j,1)$ and define 
\[\psi_j=   \frac{\chi_{B_j}}{ \sum_{k=1}^N \chi_{B_k}  }  .\]
We have $f=\sum_{j=1}^N f_j$, where $f_j=f\,\psi_j$. Since $\frac{f_j}{\mu(B_j)^{1/2}\|f_j\|_{L^2}}$ is a global atom, then 
\[
\begin{aligned}
\|f\|_{\mathfrak {h} ^1}&\leq \sum_{j=1}^N\mu(B_j)^{\frac{1}{2}}\|f_j\|_{L^2}\lesssim  \sum_{j=1}^N \|f_j\|_{L^2} \lesssim N^{\frac{1}{2}}\left(\sum_{j=1}^N    \|f_j\|_{L^2}^2\right)^{1/2}\lesssim \mu(B)^{\frac{1}{2}}\|f\|_{L^2},
\end{aligned}
\]
where we used Schwarz's inequality and the fact that $N\lesssim \mu(B)$. \hfill\ensuremath{\blacksquare}

\begin{lemma}\label{normah1funzioneL2medianulla}
Let $f$ be a function in $L^2(\mathbb X)$ with vanishing integral supported in an annulus $A_{R-r}^{R+r}$, $r\in (0,1]$, $R>r$. Then $f$ is in $\mathfrak{h}^1(\mathbb X)$ and
\[
\|f\|_{\mathfrak{h}^1}\lesssim  \log(1/r)\,e^{\rho R}r^{1/2}\,\|f\|_{L^2}.
\]
\end{lemma}
{\bf{Proof.}} We take a $r/3$-discretization $\Sigma$ of $\mathbb X$. The set $A_{R-r}^{R+r}\cap \Sigma$ has at most $N$ elements $z_1,\dots, z_N$. Then $A_{R-r}^{R+r}\subseteq \cup_{j=1}^NB_j\subseteq A_{R-2r}^{R+2r}$, so that 
\begin{equation}\label{estN}
N\leq C\,r^{-n}\mu(A_{R-2r}^{R+2r})\lesssim r^{-n+1} e^{2\rho R},
\end{equation}
the second inequality by~\eqref{misura_anello}. Let $K$ be the lowest integer such that $2^{K}r>1$ and for every $k=0, \dots, K$ and $j=1,\dots, N$ denote by $B^k_j$ the ball $B(z_j,2^kr)$ and define 
\[
\psi_j=   \frac{\chi_{B^0_j}}{ \sum_{i=1}^N \chi_{B^0_i}  },  \qquad \phi_j^k= \frac{\chi_{B^k_j}}{ \mu(B^k_j)  } .
\] 
Clearly $\int \phi^k_j \, d\mu=1$ and $\|\phi^k_j\|_{L^2} =\mu(B^k_j)^{-1/2}$. Set $f^0_j=f\psi_j$, so that $f=\sum_{j=1}^Nf^0_j$. Next, define
\[
\begin{aligned}
a_j^0
&= f_j^0-\phi^0_j  \, \int f_j^0\, d\mu,
\\
a^k_j&=(\phi^{k-1}_j-\phi^k_j)\int f_j^0\, d\mu\qquad k=1,\dots, K-1,\\
a^K_j
&= \phi^{K-1}_j  \,  \int f_j^0\, d\mu.
\end{aligned}
\]
Then, the support of $a_j^0$ is contained in $B^0_j$, the integral 
of  $a_j^0$ vanishes and  
\[
\|{a_j^0}\|_{L^2}\leq \|{f_j^0}\|_{L^2}+\mu(B_j^0)^{-1/2} \|{f_j^0}\|_{L^2}\mu(B_j^0)^{1/2} = 2 \, \|{f_j^0}\|_{L^2} \, .
\] 
Hence, by Lemma~\ref{normah1funzioneL2}
\[
\|a_j^0\|_{\mathfrak h^1} 
\leq 2  \,      \|f_j^0\|_{L^2}\mu(B^0_j)^{1/2} \, .
\]
The function $a_j^k$ is supported in $B_j^k$, the integral 
of  $a_j^k$ vanishes and  
\[
\|{a_j^k}\|_{L^2}
\leq   \|{f_j^0}\|_{L^2} \, \mu(B^0_j)^{1/2}  (\mu(B^{k-1}_j)^{-1/2} +\mu(B^{k}_j)^{-1/2})  \, .
\] 
Then, again by Lemma~\ref{normah1funzioneL2}
\[
\begin{aligned}
\|a_j^k\|_{\mathfrak h^1}&\leq  \|{f_j^0}\|_{L^2} \, \mu(B^0_j)^{1/2}\, \mu(B_j^k)^{1/2} (\mu(B^{k-1}_j)^{-1/2} +\mu(B^{k}_j)^{-1/2})\\
&= \|{f_j^0}\|_{L^2} \, \mu(B^0_j)^{1/2}\,\frac{\mu(B^{k-1}_j)^{1/2} +\mu(B^{k}_j)^{1/2}}{\mu(B^{k-1}_j)^{1/2} }\lesssim \|{f_j^0}\|_{L^2} \, \mu(B^0_j)^{1/2}.
\end{aligned}
\]
Finally, the function $a_j^K$ is supported in $B_j^K$, whose radius is bigger than $1$ but smaller than $2$, so that by Lemma~\ref{normah1funzioneL2}
\[
\begin{aligned}
\|{a_j^K}\|_{\mathfrak h^1} \lesssim \|{a_j^K}\|_{L^2}  \lesssim   \|{f_j^0}\|_{L^2} \, \mu(B^0_j)^{1/2}.
\end{aligned}
\]
It follows that $f=\sum_{j=1}^Nf^0_j=\sum_{j=1}^N\sum_{k=0}^Ka_j^k$ and
\[
\begin{aligned}
\|{f}\|_{\mathfrak h^1}  
& \lesssim \sum_{k=0}^K \sum_{j=1}^N  \|{f_j^0}\|_{L^2} \, \mu(B^0_j)^{1/2}  \\
&\leq K\,N^{1/2}    \left(\sum_{j=1}^N \|{f_j^0}\|_{L^2}^2 \right)^{1/2}r^{n/2}
\lesssim \log(1/r)\,    e^{\rho R} r^{1/2}  \|{f}\|_{L^2},
\end{aligned}
\]
the last inequality by \eqref{estN} and since $  \sum_{j=1}^N \|{f_j^0}\|_{L^2}^2\leq M   \|{f}\|_{L^2}^2 $, where $M$ is the constant in~\eqref{Muniform}. This completes the proof of the lemma.  \hfill\ensuremath{\blacksquare}

\begin{lemma}\label{aastgamma}
Let $\gamma$ be a radial function supported in $B(o,\beta)$. 
\begin{enumerate}[(i)]
\item If $a$ is a global atom at scale $1$ supported in $B(o,1)$, then 
\[
\|a\ast \gamma\|_{\mathfrak{h}^1}\leq C\, \mu(B(o,1+\beta))^{1/2}\,\|\gamma\|_{L^2}\,;
\]
\item if $a$ is a standard atom supported in $B(o,r)$, $r\in (0,1]$, then 
\[
\|a\ast \gamma\|_{\mathfrak{h}^1}\leq C\,  \mu(B(o,r+\beta))^{1/2}\,\min( \|\gamma\|_{L^2},r\,\| \nabla \gamma \|_{L^2}  ),
\]
\end{enumerate}
where $\nabla$ is the Riemannian gradient. 
\end{lemma}
{\bf{Proof.}} To prove (i), if $a$ is a  global atom supported in $B(o,1)$, then $a\ast \gamma$ is supported in $B(o,1+\beta)$ and
\[
\|a\ast\gamma\|_2 \leq \|a\|_{L^1}\,\|\gamma\|_{L^2}\leq \|\gamma\|_{L^2}.
\] 
Thus, (i) follows from Lemma~\ref{normah1funzioneL2}.

To prove (ii), if $a$ is a standard atom supported in $B(o,r)$, $r\leq 1$, then $a\ast \gamma$ is supported in $B(o,r+\beta)$ and again
\[
\|a\ast\gamma\|_2 \leq \|a\|_{L^1}\,\|\gamma\|_{L^2}\leq \|\gamma\|_{L^2}.
\] 
By arguing as in~\cite[Lemma 2.7]{MMV} and using the cancellation of the atom we obtain that
\begin{equation}\label{L2convgrad}
\|a\ast\gamma\|_2 \leq r\,\|\nabla \gamma\|_{L^2}. 
\end{equation}
Thus, (ii) follows from Lemma~\ref{normah1funzioneL2}. \hfill\ensuremath{\blacksquare}

\begin{lemma}\label{Sobolev}
Let $m$ be an even symbol in $S^b_0$ and $\mathcal U_m$ be the operator defined by the Fourier multiplier $m$. The following hold: 
\begin{enumerate}[(i)]
\item if $2\leq q<\infty$ and $\frac{1}{q}=\frac{1}{2}+\frac{b}{n}$, then $\mathcal U_m$ is bounded from $L^2(\mathbb X)$ to $L^q(\mathbb X)$;
\item if $1< s\leq 2$ and $\frac{1}{s}=\frac{1}{2}-\frac{b}{n}$, then $\mathcal U_m$ is bounded from $L^s(\mathbb X)$ to $L^2(\mathbb X)$.
\end{enumerate}
\end{lemma}
{\bf{Proof.}} Part (i) is proved in~\cite[Lemma 3]{I}. 

Part (ii) follows by a duality argument. Indeed, the adjoint of $\mathcal U_m$ is the operator $\mathcal U_{\overline{m}}$. Since $m\in S^b_0$ also $\overline{m}\in S^b_0$. By (i) the operator $\mathcal U_{\overline{m}}$ is bounded from $L^2(\mathbb X)$ to $L^q(\mathbb X)$, with $2\leq q<\infty$ and $\frac{1}{q}=\frac{1}{2}+\frac{b}{n}$. Then $\mathcal U_m$ is bounded 
from $L^{q'}(\mathbb X)$ to $L^2(\mathbb X)$. Let $s=q'$. Then $1< s\leq 2$ and $\frac{1}{s}=1-\frac{1}{q}=1-\frac{1}{2}-\frac{b}{n}=\frac{1}{2}-\frac{b}{n}$, as required.  \hfill\ensuremath{\blacksquare}

\section{Boundedness of $\mathcal T_{t}$ from $\mathfrak{h}^1(\mathbb X)$ to $L^1(\mathbb{X})$} \label{main1}
In this section, we prove part (i) of Theorem~\ref{main_teo_intro}. The proof is inspired to that of~\cite[Proposition 4]{I}.

\smallskip 

\noindent {\bf{Proof of Theorem~\ref{main_teo_intro}~(i).}} By Proposition~\ref{uniform} and since $\mathcal T_t$ is left invariant it is enough to prove that 
\[
\sup\{ \|\mathcal T_t a\|_{L^1}\colon a \mbox{ $\mathfrak{h}^1$-atom supported in $B(o,r)$, } r\leq 1\}\lesssim e^{\rho \,t}.
\]
Let $a$ be an atom supported in $B(o,r)$, $r\leq 1$. We separate two different cases, according to the values of $t$.

\smallskip

{\bf{{\emph{Case I: $t\geq 1/2$.}}}}  We define the set
\[
B^*\coloneqq \{ x\in \mathbb{X} \colon | |x|-t|<10 r \},
\]
whose measure is $\mu(B^*)\lesssim r e^{2\rho t}$, and split
\[
\| \mathcal{T}_t a\|_{L^1} = \| \mathcal{T}_t a\|_{L^1(B^*)} + \| \mathcal{T}_t a\|_{L^1((B^*)^c)}.
\]
We observe that by H\"older inequality
\[
\| \mathcal{T}_t a\|_{L^1(B^*)} \leq \mu(B^*)^{1/2}\| \mathcal{T}_t a \|_{L^2}\lesssim e^{\rho t} r^{1/2}\| \mathcal{T}_t a \|_{L^2}.
\]
Moreover, by Lemma~\ref{Sobolev} (ii) with $\frac{1}{s}=\frac12-\left(-\frac{d}{n}\right)=\frac12+\frac{n-1}{2n}=1-\frac{1}{2n}$, H\"older inequality and the size condition of the atom
\begin{equation}\label{normaL2TtaL1}
\|\mathcal T_ta\|_{L^2} \lesssim \|a\|_{L^s} \lesssim \mu(B)^{-1+1/s} \lesssim r^{-1/2}.
\end{equation}
Thus $\| \mathcal{T}_t a\|_{L^1(B^*)} \lesssim e^{\rho t}$. 
 
Let now $k_{t}$ be the radial kernel of the operator $\mathcal T_{t}$, and let $K_{t}$ be the function on $[0,\infty)$ such that $k_{ t}(x)=K_{ t}(|x|)$. It remains to estimate the $L^1$-norm of $a* k_t$ on $(B^*)^c$. In order to do this, we take a function
\[
\psi_t \in C_c^\infty(\mathbb{X}), \quad \psi_t(x)=1 \; \mbox{ if }\; ||x|-t| <\tfrac{1}{10}, \quad \psi_t(x)=0\; \mbox{ if }\; ||x|-t| \geq \tfrac{2}{10},
\]
with values in $[0,1]$, define $\Psi_t(|x|) = \psi_t(x)$, and split the kernel $k_t$ in its singular part $s_t$ and its good part $g_t$ as
\[
k_t = k_t \psi_t + k_t(1-\psi_t) \eqqcolon s_t + g_t.
\]
Observe that this induces a splitting $K_t = K_t\Psi_t + K_t(1-\Psi_t)\eqqcolon S_t + G_t$ of functions defined on $\mathbb{R}^+$. It is proved in~\cite[p.\ 287]{I} that 
\begin{equation}\label{GtIonescu}
|G_t(s)|\lesssim 
\begin{cases}
s^{-d-1}& \mbox{if } s\leq \tfrac{1}{10}\\
e^{-\rho s}|t-s|^{-2}& \mbox{if } \tfrac{1}{10} \leq s \leq t-\frac{1}{10}\\
e^{\rho t} e^{-2\rho s}|t-s|^{-2} & \mbox{if } s\geq t+\frac{1}{10}\\
\end{cases}
\end{equation}
from which $\|g_t\|_{L^1}\leq e^{\rho t}$. Thus
\[
\| a * g_t\|_{L^1((B^*)^c)}\leq \| a * g_t\|_{L^1} \leq \| a\|_{L^1} \|g_t\|_{L^1}\leq  e^{\rho t}.
\]
As for the convolution with $s_t$, we first consider the case when $a$ is a global atom. Since $\psi_t$ is supported in the annulus $A_{t-2/10}^{t+2/10}$, the convolution $a*s_t$ is supported in the annulus $A_{t-6/5}^{t+6/5}$. Then by H\"older inequality
\[
\|a*s_t\|_{L^1((B^*)^c)}\leq \|a*s_t\|_{L^1}\lesssim \mu(A_{t-6/5}^{t+6/5})^{1/2}\|a*s_t\|_{L^2} \lesssim e^{\rho t}\|a*s_t\|_{L^2} 
\]
where
\[
\|a*s_t\|_{L^2} \lesssim \|a*k_t\|_{L^2} + \| a * g_t\|_{L^2} \lesssim \|\mathcal{T}_t\|_{L^2\rightarrow L^2}\| a\|_{L^2} +\| g_t\|_{L^2} \|a\|_{L^1} \lesssim 1,
\]
since $\|g_t\|_{L^2}\lesssim 1$ by~\eqref{GtIonescu}. Thus $\|a*s_t\|_{L^1((B^*)^c)}\lesssim e^{\rho t}$. If instead $a$ is a standard atom, by its cancellation condition it is easy to see that
\[
a\ast s_t(x) = \int_G a(z) \left[ s_t(z^{-1}x) - s_t(x)\right]\, dz = \int_B a(z) \left[ S_t(|z^{-1}x|) - S_t(|x|)\right]\, dz
\] 
for every $x\in \mathbb{X}$, so that
\[
\| a *s_t\|_{L^1((B^*)^c)}\leq \int_B |a(z)| \int_{(B^*)^c}  |S_t(|z^{-1}x|)-S_t(|x|)|\, dx \, dz.
\]
It remains to observe that, since $\left| \partial_s S_t (s) \right|\lesssim e^{- \rho t} |t-s|^{-2}$ as shown in~\cite[p.\ 287]{I},
\begin{align*}
\sup_{z\in B}\;  \int_{(B^*)^c}  |S_t(|z^{-1}x|)-S_t(|x|)|\, dx
&\lesssim \sup_{z\in B}\; |z| \int_{10 r\leq ||x|-t| \leq r+2/10} \left| \partial_s S_t (|x|) \right| \, dx
\\& \lesssim r e^{-\rho t} \int_{10 r\leq ||x|-t| \leq r+2/10} ||x|-t|^{-2} \, dx  \lesssim e^{\rho t},
\end{align*}
which concludes the proof of the Case I.
\smallskip

{\bf{{\emph{Case II: $t< 1/2$.}}}}  After defining the set
\[
B^*\coloneqq \{ x\in \mathbb{X} \colon | |x|-t|<10 r \} \cup B(0,10r),
\]
we proceed as in the previous case. Since $\mu(B^*)\lesssim r$, we get $\|\mathcal{T}_t a\|_{L^1(B^*)} \leq C $ again by~\eqref{normaL2TtaL1}. In order to estimate $\|\mathcal{T}_t a\|_{L^1((B^*)^c)} $, we pick a function
\[
\psi_0 \in C_c^\infty(\mathbb{X}), \quad \psi_0(x)=1 \; \mbox{ if }\; |x|\leq \tfrac{3}{4}, \quad \psi_0(x)=0\; \mbox{ if }\; |x|\geq 1,
\]
and split again the kernel $k_t$ as
\[
k_t = k_t \psi_0 + k_t(1-\psi_0) = s_t + g_t.
\]
We let $\Psi_0$, $S_t$ and $G_t$ be the associated functions on $\mathbb{R}^+$. It is proved in~\cite[p.\ 288]{I} that 
\[
|G_t(s)|\lesssim e^{-2\rho s}|t-s|^{-2} \qquad \forall \, s\geq \tfrac{3}{4},
\]
so that $\|g_t\|_{L^1}\leq C$, hence $\| a * g_t\|_{L^1((B^*)^c)} \leq C$. As for the convolution with $s_t$, if $a$ is a global atom then as before
\[
\|a*s_t\|_{L^1}\lesssim \mu(A_{t-1}^{t+1})^{1/2}\|a*s_t\|_{L^2} \lesssim \|\mathcal{T}_t\|_{L^2\rightarrow L^2}\| a\|_{L^2} + \| g_t\|_{L^2} \|a\|_{L^1}\lesssim 1,
\]
while if $a$ is a standard atom, by its cancellation condition we obtain again
\[
\| a *s_t\|_{L^1((B^*)^c)}\leq \int_B |a(z)| \int_{(B^*)^c}  |S_t(|z^{-1}x|)-S_t(|x|)|\, dx \, dz.
\]
Proceeding as in~\cite[p.\ 288]{I}, $S_t$ may be written as the sum of two functions $S_{1,t} +S_{2,t}$ such that $S_{1,t}(s)\leq s^{-d-1}$ (hence $s_{1,t}\in L^1(\mathbb{X})$) while
\[
\left| \partial_s S_{2,t} (s) \right|\lesssim s^{-d }\left(|t-s|^{-2} + s|t-s|^{-1}\right).
\]
The proof may be completed as before. \hfill\ensuremath{\blacksquare}

\section{Boundedness of $\mathcal T_{t}$ on $\mathfrak{h}^1(\mathbb X)$}\label{main}
In this section we prove part (ii) of Theorem~\ref{main_teo_intro}, but first we need some preliminary results. We recall the behavior of the Harish-Chandra function and of spherical functions on noncompact symmetric spaces of rank one. It follows from~\cite[Propositions A.1, A.2]{I} and is based on various results in~\cite{ST}. We denote by $\rho'$ the number $\rho+\frac{1}{10}$. 
\begin{lemma}\label{sphericalfunctions}
The Harish-Chandra function ${\bf{c}}$ satisfies the following:
\begin{enumerate}[(i)]
\item for all $\lambda\in\mathbb R$ 
\[
|{\bf{c}}(\lambda)|^{-2}={\bf{c}}(\lambda)^{-1}\,{\bf{c}}(-\lambda)^{-1}\,;
\]
\item the function $\lambda\mapsto \lambda^{-1}\,{\bf{c}}(-\lambda)^{-1}$ is analytic inside the region $   \mathrm{Im}\, \lambda \geq 0$ and for all $\alpha \geq 0$ there exists a positive constant $C_{\alpha}$ such that 
\[
\Big|    \partial_{\lambda}^{\alpha} (\lambda^{-1}\,{\bf{c}}(-\lambda)^{-1})\Big|\leq C_{\alpha}\,(1+|\mathrm{Re}\,\lambda|)^{d-1-\alpha}\qquad \forall  \; 0\leq \mathrm{Im}\, \lambda \leq \rho';
\]
\item the function $\lambda\mapsto \lambda\,{\bf{c}}(\lambda)$ is analytic in a neighborhood of the real axis and for all $\alpha \geq 0$ there exists a positive constant $C_{\alpha}$ such that 
\[
\Big|    \partial_{\lambda}^{\alpha}  (\lambda\,{\bf{c}}(\lambda))\Big|\leq C_{\alpha}\,(1+|\mathrm{Re}\,\lambda|)^{1-d-\alpha}\qquad \forall \; \lambda\in\mathbb R.\]
\end{enumerate}
The spherical functions $\phi_{\lambda}$ satisfy the following properties:
\begin{enumerate}[(a)]
\item $| {\partial_s^{\ell}}\phi_{\lambda}(s)|\leq C\,e^{-\rho s}(1+s)\,(1+|\lambda|)^{\ell}\qquad \forall \; \lambda,s\in\mathbb R$, $\ell\in\mathbb N$.
\item If $s\leq 1$, $\lambda\in\mathbb R$ and $s|\lambda|\geq 1$, for every $N\in\mathbb N$, $\phi_{\lambda}$ can be written as
\[
\phi_{\lambda}(s)=e^{i\lambda s}a_1(\lambda,s)+e^{-i\lambda s}a_1(-\lambda,s)+O(\lambda,s),
\] 
where the functions $a_1,O:\{(s,\lambda)\in \mathbb R\times [0,1] \colon s|\lambda|\geq 1\}\rightarrow \mathbb C$ satisfy
\[
\Big| {\partial_{\lambda}^{\alpha}} {\partial_s^{\ell}}a_1(\lambda,s)\Big|\leq C\,[s(1+|\lambda|)]^{-d}\,s^{-\ell}\,(1+|\lambda|)^{-\alpha}\qquad \ell\in \{0,1\}, \, \alpha\in [0,N]
\]
and
\[
| {\partial_s^{\ell}}O(\lambda,s)|\leq C\,[s(1+|\lambda|)]^{-d-N-1-\ell}.
\]
\item If $s\geq 1/10$, then 
\[
\phi_{\lambda}(s)=e^{-\rho s}\left(e^{i\lambda s}{\bf{c}}(\lambda)a_2(\lambda,s)+e^{-i\lambda s}{\bf{c}}(-\lambda)a_2(-\lambda,s)\right),
\]
where the function $a_2$ is such that for all $\alpha\geq 0$ there exist positive constants $C_{\alpha }$ such that 
\[
\Big|  {\partial_{\lambda}^{\alpha}} {\partial_s^{\ell}}a_2(\lambda,s)\Big|\leq C_{\alpha }\,(1+|\mathrm{Re}\,\lambda|)^{-\alpha}\qquad \forall \, \ell\in\{0,1\},\; s\geq \tfrac{1}{10}, \; 0\leq \mathrm{Im}\, \lambda \leq  \rho'.
\]
\end{enumerate}
\end{lemma}
{\bf{Proof.}} The properties of the Harish--Chandra function were given in~\cite{I}. See also~\cite[Formula (2.2.5)]{AJi}. 

Formula (a) follows from~\cite[Formula 5.1.18]{GV}.

The proof of (b) follows the same outline of the proof of~\cite[Proposition A.2 (b)]{I}. The only difference is that following the same arguments it is possible to estimate the derivatives of the term $O(\lambda,s)$ which were not estimated in~\cite{I}. 

The proof of (c) is given in~\cite[Proposition A.2 (c)]{I}.  \hfill\ensuremath{\blacksquare} 
 
In the following proposition we shall prove pointwise estimates of the kernel of the operator $\mathcal T_{t}$ and of its derivative. We will distinguish the cases when $t$ is either large or small. Let us mention that Ionescu~\cite{I} estimated the kernels of the operator $\mathcal T_{t}$ (but not their derivatives) far from the sphere of radius $t$, while he gave estimates of the derivatives of the kernels (but not of the kernels) near the sphere of radius $t$.

\begin{proposition}\label{prop_kt}
Let $\epsilon >0$ and $m\in S^{-d-\epsilon}_{\rho}$ be an even symbol. Let $k_{t}$ be the radial kernel of the operator $\mathcal T_{t}$ and $K_{t}$ be the function on $[0,\infty)$ such that $k_{ t}(x)=K_{ t}(|x|)$. 

\noindent If $t\geq \frac12$, then 
\begin{equation}\label{Kttgrande}
\begin{aligned}
|K_{ t}(s)|&\lesssim\begin{cases}
s^{-d-1+\epsilon}&s\leq \frac{1}{10}\\
e^{-\rho s} |t-s|^{-2+[\epsilon]} & \frac{1}{10}\leq s\leq t-\frac{2}{10}\\
e^{-\rho t} |t-s|^{-1+\epsilon} &    t-\frac{2}{10}\leq s\leq t+\frac{2}{10}     \\
e^{\rho t}\,e^{-2\rho s}\,|t-s|^{-2+[\epsilon]}  &s\geq t+\frac{2}{10};
\end{cases}
\end{aligned}
\end{equation}
\begin{equation}\label{Kt'tgrande}
 \begin{aligned}
|K_{t}'(s)|&\lesssim\begin{cases}
s^{-d-2+\epsilon}&s\leq \frac{1}{10}\\
e^{-\rho s} |t-s|^{-2+[\epsilon]}& \frac{1}{10}\leq s\leq t-\frac{2}{10}\\
e^{-\rho t} |t-s|^{-2+\epsilon} &    t-\frac{2}{10}\leq s\leq t+\frac{2}{10}    \\
 e^{\rho t}\,e^{-2\rho s}\,|t-s|^{-2+[\epsilon]} &s\geq t+\frac{2}{10}.
\end{cases}
\end{aligned}
\end{equation}
If $t< \frac{1}{2}$, then 
\begin{equation}\label{Kttpiccolo}
 \begin{aligned}
 |K_{ t}(s)|&\lesssim \begin{cases}
 e^{-2\rho s}\,|t-s|^{-2+[\epsilon]}  &s\geq 1\\
 s^{-d-1+\epsilon}+s^{-d}|t-s|^{-1+\epsilon}&s\leq 1\,;
 \end{cases}
 \end{aligned}
 \end{equation}
\begin{equation}\label{Kt'tpiccolo}
 \begin{aligned}
  |K_{ t}'(s)|&\lesssim \begin{cases}
 e^{-2\rho s}\,|t-s|^{-2+[\epsilon]} &s\geq 1\\
  s^{-d-2+\epsilon}+s^{-d}|t-s|^{-2+\epsilon}+s^{-d-1}|t-s|^{-1+\epsilon}&s\leq 1.
 \end{cases} 
\end{aligned}
\end{equation}
\end{proposition}

{\bf{Proof.}} Since the operator $\mathcal T_{t}$ corresponds to the spherical Fourier multiplier $\lambda\mapsto m(\lambda)\cos({t \lambda })$, by the inversion formula for the spherical transform~\eqref{inversion} we get

\begin{equation}\label{kt}
K_{t}(s)=C\,\int_{\mathbb R}m(\lambda)\,\cos(t \lambda)\,\phi_{\lambda}(s)\,|{\bf{c}}(\lambda)|^{-2}\,d\lambda.
\end{equation}
We distinguish the cases when $t$ is either large or small.
 
\smallskip 
 
{\it{\textbf{Case I: $t\geq 1/2$}.}}  Let $\Psi_t$ be a smooth cutoff function such that 
\[
\Psi_t(s)=1 \; \mbox{ if }\; |s-t|\leq \tfrac{1}{10}, \quad \Psi_t(s)=0\; \mbox{ if }\; |s-t|\geq \tfrac{2}{10}.
\]
Let $S_t\coloneqq \Psi_t\,K_{t}$ and $G_t\coloneqq(1-\Psi_t)\,K_{t}$. To prove~\eqref{Kttgrande} and~\eqref{Kt'tgrande} it is enough to estimate $S_t$ and $G_t$ and their derivatives. We shall repeatedly use, without further mention,~\cite[Lemma A.2]{APV1} to estimate the Fourier transform of a symbol of some given order.

\smallskip

We first consider $S_t$. Observe that $S_t(s)= 0$ unless $|t-s|\leq \frac{2}{10}$, i.e.\ $t-\tfrac{2}{10}\leq s\leq t+\tfrac{2}{10}$. 
From~\eqref{kt} and Lemma~\ref{sphericalfunctions} we deduce that 
\[
S_t(s)= C\,\Psi_t(s)\,e^{-\rho\,s }\,\int_{\mathbb R}m(\lambda)\,\cos({t \lambda})\,e^{i\lambda s}\,a_2(\lambda,s)\,{\bf{c}}(-\lambda)^{-1}\,d\lambda.
\]
Since by Lemma~\ref{sphericalfunctions} ({c}) the function $\lambda \mapsto m(\lambda)\,\,a_2(\lambda,s)\,{\bf{c}}(-\lambda)^{-1}$ is a symbol on the real line of order $-\epsilon$
\[
|S_t(s)|\lesssim e^{-\rho t}\, |t-s|^{-1+\epsilon} .
\]
Similarly, one can see that $|S_t'(s)|\lesssim e^{-\rho t}\, |t-s|^{-2+\epsilon}$. 

\smallskip

To estimate $G_t$ and its derivative, we observe that $G_t(s)=0$ unless $|t-s|\geq \tfrac{1}{10}$. The function $G_t$ can be estimated as in~\cite[Formula (3.9)]{I} (see also~\eqref{GtIonescu}). To estimate the derivative of $G_t$ we distinguish different cases. 

We first consider the case when $s\leq \frac{1}{10}$. We choose a smooth cutoff function $\eta$ such that 
\[
\eta(v)=1 \; \mbox{ if } \; |v|\leq 1, \quad \eta(v)=0 \; \mbox{ if } \;  |v|\geq 2.
\]
By Lemma~\ref{sphericalfunctions} (b) we write
\[
\begin{aligned}
G_t(s)&=C\,\left(1-\Psi_t(s)\right)\,\int_{\mathbb R}  \eta(\lambda\, s)\, \phi_{\lambda}(s)\, m(\lambda)\,\cos({t\lambda })\,|{\bf{c}}(\lambda)|^{-2}\,d\lambda\\
&   +C\,\left(1-\Psi_t(s)\right)\,\int_{\mathbb R}\,  \left(1-\eta(\lambda s)\right)\,O(\lambda,s)\, m(\lambda)\,\cos({t\lambda }) \,|{\bf{c}}(\lambda)|^{-2}\,d\lambda\\
&+C\,\left(1-\Psi_t(s)\right)\,\int_{\mathbb R}\left(1-\eta(\lambda s)\right)e^{i\lambda s}\,a_1(\lambda, s)\, m(\lambda)\,\cos({t\lambda })\,|{\bf{c}}(\lambda)|^{-2}\,d\lambda.
\end{aligned}
\]
Then
\begin{equation}\label{Btlocale}
\begin{aligned}
G_t'(s)&=C\, \int_{\mathbb R}  \Big[- \Psi'_t(s)\eta(\lambda\, s)\, \phi_{\lambda}(s)\, +\left(1-\Psi_t(s)\right)\lambda \eta'(\lambda s)\,\phi_{\lambda}(s)\\
&\quad +\left(1-\Psi_t(s)\right)  \eta(\lambda s)\, {\partial_s}\phi_{\lambda}(s) \Big]
m(\lambda)\, \cos({t\lambda })\,|{\bf{c}}(\lambda)|^{-2}\,d\lambda\\
&   \quad +C\, \int_{\mathbb R}\,   \Big[- \Psi'_t(s)\left(1-\eta(\lambda\, s)\right)\, O({\lambda},s)\, -\left(1-\Psi_t(s)\right)\lambda \eta'(\lambda s)\,O(\lambda,s)\\
&\quad +\left(1-\Psi_t(s)\right) \left(1- \eta(\lambda s)\right)\, {\partial_s}O(\lambda,s) \Big]
 \, m(\lambda)\, \cos({t\lambda })\,|{\bf{c}}(\lambda)|^{-2}\,d\lambda\\
&\quad +C\, \,\int_{\mathbb R}   \Big[-\Psi'_t(s)\left(1-\eta(\lambda\, s)\right)\, a_1(\lambda, s)\, \, -\left(1-\Psi_t(s)\right)\lambda \eta'(\lambda s)\,  a_1(\lambda, s)\,\\
&\quad +\left(1-\Psi_t(s)\right) \left(1- \eta(\lambda s)\right)\, i\lambda\, a_1(\lambda, s)\, +\left(1-\Psi_t(s)\right) \left(1- \eta(\lambda s)\right)\,   {\partial_s}a_1(\lambda, s)\, 
\Big]\\
&\quad \times e^{i\lambda s}\,m(\lambda)\, \cos({t\lambda })\,
\,|{\bf{c}}(\lambda)|^{-2}\,d\lambda\\
&=G'_{1,t}(s)+G'_{2,t}(s)+G'_{3,t}(s).
\end{aligned}
\end{equation}
By Lemma~\ref{sphericalfunctions} (a)  
\[
\begin{aligned}
|G'_{1,t}(s)|&\lesssim \int_0^{2/s}(1+\lambda)^{d-\epsilon +1}\, d\lambda \lesssim s^{-d-2+\epsilon}.
\end{aligned}
\]
Similarly, by Lemma~\ref{sphericalfunctions} (b), (with $N=0$)
\[
\begin{aligned}
|G'_{2,t}(s)|&\lesssim  \int_{1/s}^{\infty}\,  s^{-d-1}\,\lambda^{-\epsilon-1}\,d\lambda  \lesssim s^{-d-1+\epsilon}.
\end{aligned}
\]
To estimate $G'_{3,t}$ we write $\cos(t\lambda)=(e^{it\lambda}+e^{-it\lambda})/2$ and integrate by parts twice:
\[
\begin{aligned}
|G'_{3,t}(s)|&\lesssim \frac{1}{|t-s|^2}    \int_{\mathbb R} \Big| \partial^2_{\lambda} \Big[-\Psi'_t(s)\left(1-\eta(\lambda\, s)\right)\, a_1(\lambda, s)\, \, -\left(1-\Psi_t(s)\right)\lambda \eta'(\lambda s)\,  a_1(\lambda, s)\,\\
&+\left(1-\Psi_t(s)\right) \left(1- \eta(\lambda s)\right)\, i\lambda\, a_1(\lambda, s)\, +\left(1-\Psi_t(s)\right) \left(1- \eta(\lambda s)\right)\,   {\partial_s}a_1(\lambda, s)\, 
\Big]\\
&\times  \,m(\lambda)\, \,|{\bf{c}}(\lambda)|^{-2}\, \Big| \, d\lambda.
\end{aligned}
\]
By applying Lemma~\ref{sphericalfunctions} (b), we can easily show that $|G'_{3,t}(s)|\lesssim s^{-d-2+\epsilon}$. 

Thus from~\eqref{Btlocale} and the estimates above, we deduce that for every $s\leq \frac{1}{10}$, $|G_t'(s)|\lesssim s^{-d-2+\epsilon}. $ 

We now consider the case $s\geq \frac{1}{10}$. By Lemma~\ref{sphericalfunctions} ({c}) we have
\begin{equation}\label{Bt}
G_t(s)=C\,(1-\Psi_t(s))\,e^{-\rho s}\,\int_{\mathbb R} m(\lambda)\,a_2(\lambda,s)\,{\bf{c}}(-\lambda)^{-1}\,e^{i\lambda s}\, \cos({t\lambda })\,d\lambda,
\end{equation}
so that
\[
\begin{aligned}
G_t'(s)&=-C  \Psi'_t(s)\,\,e^{-\rho s}\,\int_{\mathbb R} m(\lambda)\,a_2(\lambda,s)\,{\bf{c}}(-\lambda)^{-1}\,e^{i\lambda s}\,\cos({t\lambda })\,d\lambda\\
&\quad+C\,(1-\Psi_t(s))\,e^{-\rho s}\,\int_{\mathbb R} m(\lambda)\,\left( {\partial_s}a_2(\lambda,s)+(-\rho+i\lambda) \,a_2(\lambda,s)\right)\\&\quad \times {\bf{c}}(-\lambda)^{-1}\,e^{i\lambda s}\,\cos({t\lambda })\,d\lambda.
\end{aligned}
\]
Since $\lambda\mapsto m(\lambda)\,\left(a_2(\lambda,s)+ {\partial_s}a_2(\lambda,s)\right)\,{\bf{c}}(-\lambda)^{-1}$ is a symbol of order $-\epsilon$, and $\lambda\mapsto m(\lambda)\,i\lambda \,a_2(\lambda,s)\,{\bf{c}}(-\lambda)^{-1}$ is a symbol of order $1-\epsilon$, we obtain that 
\[
|G_t'(s)|\lesssim  e^{-\rho s} |t-s|^{-2+[\epsilon]} \qquad \mbox{ if }\; \tfrac{1}{10}\leq s\leq t-\tfrac{1}{10}. 
\]
It remains to consider the case when $s\geq t+\frac{1}{10}$. In order to do this, we move the contour of integration in formula~\eqref{Bt} to the line $\mathbb R+i\rho $ and obtain 
\[
\begin{aligned}
G_t(s)&=\,C\,(1-\Psi_t(s))\,e^{-2\rho s}\,\\ &\quad \times  \int_{\mathbb R} m(\lambda+i\rho)\,a_2(\lambda+i\rho,s)\,{\bf{c}}(-\lambda-i\rho)^{-1}\,e^{i\lambda s} \cos(t{(\lambda +i\rho)})\,d\lambda.
\end{aligned}
\]
By taking the derivative we get
\[
\begin{aligned}
G_t'(s)&= e^{-2 \rho s} \int_\mathbb{R} m(\lambda +i\rho){\bf{c}}(-\lambda-i\rho)^{-1} \cos(t(\lambda +i\rho)) e^{i\lambda s} \Big\{ -\Psi_t'(s)a_2(\lambda+i\rho, s) \\ 
&\qquad \qquad + (1-\Psi_t(s)) \left[ a_2(\lambda+i\rho,s )(-2\rho +i\lambda) + \partial_s a_2(\lambda +i\rho,s) \right] \Big\} \, d\lambda.
\end{aligned}
\]
The estimates of the derivatives of $a_2$ and ${\bf{c}}^{-1}$ contained in Lemma~\ref{sphericalfunctions} imply that 
\[
| G_t'(s) |\lesssim e^{-2\rho s}\,e^{\rho  t}\,|t-s|^{-2+[\epsilon]}.
\]
By combining the estimates of $G_t, G_t', S_t$ and $S_t'$ one deduces the required estimates of $K_{t}$ and its first derivative for $t$ large. 

\smallskip

{\it{\textbf{Case II: $t<1/2$.}}} Let $\Psi_0$ be a smooth cutoff function such that 
\[
\Psi_0({s})=1\; \mbox{ if } \; s\leq \tfrac{3}{4}, \qquad \Psi_0({s})=0 \; \mbox{ if } \; s\geq 1.
\]
Let $S_t=\Psi_0\,K_t$ and $G_t=(1-\Psi_0)\,K_t$. 

\smallskip

We first analyse $G_t$ and notice that $G_t(s)=0$ if $s\leq 3/4$. If $s>3/4$, then by Lemma~\ref{sphericalfunctions} ({c})
\[
G_t(s)=C\,\left(1-\Psi_0(s)\right)\,\int_{\mathbb R}m(\lambda)\,\cos(t \lambda)e^{-\rho s}\,e^{i\lambda s}\,a_2(\lambda,s)\,{\bf{c}}(-\lambda)^{-1}d\lambda,
\] 
which by moving the contour of integration from the real line to $\mathbb R+i\rho $ becomes
\[
\begin{aligned}
G_t(s)&=\,C\,\left(1-\Psi_0(s)\right)\,e^{-2\rho s}\\&\quad\times \int_{\mathbb R}m(\lambda+i\rho )\,e^{i\lambda s}\,a_2(\lambda+i\rho ,s)\,
 \cos({t(\lambda+i\rho )})\,{\bf{c}}(-\lambda-i\rho )^{-1}d\lambda.
\end{aligned}
\]
The function $G_t$ can be estimated as in~\cite[p. 289]{I}. Since $G_t$ is the Fourier transform at $s\pm t$ of a symbol of order $-\epsilon$, $s> 3/4$ and $t<1/2$,
\[
|  G_t'(s) |\lesssim  e^{-2\rho s}\,|t-s|^{-2+[\epsilon]}.
\]
\smallskip
It remains to consider $S_t$. Observe that $S_t(s)=0$ unless $s\leq 1$, hence we use Lemma~\ref{sphericalfunctions} ({c}) (with $N=0$) to write
\[
\begin{aligned}
S_t(s)&=\Psi_0(s)\int\eta(\lambda s)\,\phi_{\lambda}(s)\,m(\lambda)\, \cos(t{\lambda })\,|{\bf{c}}(\lambda)|^{-2}\,d\lambda\\
&\quad+\Psi_0(s)\int (1-\eta(\lambda s))\,O(\lambda,s)\, m(\lambda)\,  \cos(t{\lambda })\,|{\bf{c}}(\lambda)|^{-2}\,d\lambda\\
&\quad +\Psi_0(s)\int (1-\eta(\lambda s))e^{i\lambda s}\,a_1(\lambda,s)\, m(\lambda)\,  \cos(t{\lambda })\,|{\bf{c}}(\lambda)|^{-2}\,d\lambda\\
&=S_{1,t}(s)+S_{2,t}(s)+S_{3,t}(s),
\end{aligned}
\]
where $\eta$ is a smooth cutoff function such that $\eta(v)=1$ if $|v|\leq 1$ and $\eta(v)=0$ if $|v|\geq 2$. For every $s\leq 1$ we have
\[
\begin{aligned}
|S_{1,t}(s)| \lesssim \int_0^{2/s} \lambda^{d -\epsilon}d\lambda \lesssim s^{-d -1+\epsilon}
\end{aligned}
\]
and
\[
\begin{aligned}
| S_{2,t}(s)|\lesssim \int_{1/s}^{\infty}  (s\lambda)^{-d-1}  \lambda^{d -\epsilon }d\lambda \lesssim s^{-d -1+\epsilon}. 
\end{aligned}
\]
Finally, $S_{3,t}$ is the inverse Fourier transform computed at $s\pm t$ of the symbol $\lambda\mapsto (1-\eta(\lambda s)) \,a_1(\lambda,s)\, m(\lambda)\,|{\bf{c}}(\lambda)|^{-2}$ of order $-\epsilon$. Then  
\[
|  S_{3,t}(s) |\lesssim s^{-d}\,|t-s|^{-1+\epsilon}.
\]
It then follows that for every $s\leq 1$ 
\[
|S_t(s)|\lesssim 
s^{-d-1+\epsilon}+s^{-d}|t-s|^{-1+\epsilon}  .
\]
In a similar way, one can prove that
for every $s\leq 1$ 
\[
|S_t'(s)|\lesssim 
s^{-d-2+\epsilon}+s^{-d}|t-s|^{-2+\epsilon} +s^{-d-1}|t-s|^{-1+\epsilon}.
\]
By combining the estimates of $G_t, G_t', S_t$ and $S_t'$ one deduces the required estimates of $K_{t}$ and its first derivative for $t$ small.  \hfill\ensuremath{\blacksquare}

\begin{remark}
Let us notice that the kernel $K_t$ and its derivative behave in the same way far from the singularities, i.e.\ far from the point $o$ and the sphere of radius $t$, while they have a different behavior near $o$ and near the sphere of radius $t$. Observe moreover that, when $t\geq \tfrac{1}{2}$ and either $\frac{1}{10}\leq s\leq t-\frac{2}{10}$ or $s\geq t+\tfrac{2}{10}$, or when $t< \tfrac{1}{2}$ and $s\geq 1$, the power $|t-s|^{-2+[\epsilon]}$ in the estimates of $K_t(s)$ and $K'_t(s)$ may be replaced with $|t-s|^{-M}$ for any integer $M\geq -2+[\epsilon]$, provided the constant (which might depend on $M$) is properly chosen. This is a consequence of~\cite[Lemma A.2]{APV1}.
\end{remark}

We are now in the position to prove the part (ii) of Theorem~\ref{main_teo_intro}. The strategy we shall adopt consists in decomposing the kernel $k_t$ of $\mathcal{T}_t$ into a sum of compactly supported functions which we shall consider separately. We also treat separately the cases when $t$ is either large or small. The proof turns out to be more delicate when $a$ is a standard atom supported in a ball of small radius, and in this case the cancellation condition of the atom is crucial together with the estimates of the derivative of the kernel. When the atom is either a global atom or a standard atom supported in a ball of radius not too small when compared with $t$ and $1$, instead, the cancellation  of the atom plays no role and only the estimates of the kernel are involved.

\smallskip

In order to do this, we shall repeatedly use smooth cutoff radial functions, which are introduced below. We fix $r\in (0,1]$ and $t>0$.

Take a function $\phi\in C^{\infty}_c(\mathbb R)$ supported in $[1/2,2]$ such that $0\leq \phi\leq 1$, $\phi=1$ in $[1,3/2]$, $\phi(s)=1-\phi(s/2)$ for every $s\in (1,2)$ and $|\phi'|\leq C$. For every $i\in\mathbb N$ and every $x\in\mathbb X$ define 
\begin{equation}\label{phii}
\phi_i(x)=\phi\left(\frac{|x|}{2^ir}\right).
\end{equation}
Observe that $\phi_i$ is supported in the annulus $A_{2^{i-1}r}^{2^{i+1}r}$, $0\leq \phi_i\leq 1$ and $|\nabla \phi_i|\leq C\,(2^ir)^{-1}$. 

For every $h\in\mathbb N$ and $x\in\mathbb X$ define 
\begin{equation}\label{etahomegah}
\eta_h(x)= \phi\left(\frac{t-|x|}{2^h r}\right), \qquad \omega_h(x)=\phi\left(\frac{|x|-t}{2^h r}\right).
\end{equation}
The function $\eta_h$ is supported in $A^{t-2^{h-1}r}_{t-2^{h+1}r}$, $0\leq \eta_h\leq 1$ and $ |\nabla \eta_h|\leq C\,(2^h r)^{-1}$. Similarly, $\omega_h$ is supported in $A^{t+2^{h+1}r}_{t+2^{h-1}r}$, $0\leq \omega_h\leq 1$ and $|  \nabla \omega_h|\leq C\,(2^h r)^{-1}$. 

Finally, take a function $\psi\in C^{\infty}_c(\mathbb R)$ supported in $[0,2]$ such that $0\leq \psi\leq 1$, $\psi=1$ in $[2/3,4/3]$ and $\psi(s+1)=1-\psi(s)$ for every $s\in (0,1)$. For every $j\geq 2$ and $x\in\mathbb X$ define
\begin{equation}\label{psij}
\psi_j(x)=\psi\left(|x|-j+1\right).
\end{equation} 
The function $\psi_j$ is supported in $A_{j-1}^{j+1} $ and $0\leq \psi_j\leq 1$. 

\smallskip

{\bf{Proof of Theorem~\ref{main_teo_intro}~(ii)}}. By Proposition~\ref{uniform} and the left invariance of $\mathcal T_t$ it is enough to prove that 
\[
\sup\{ \|\mathcal T_t a\|_{\mathfrak h^1}: a \mbox{ $\mathfrak{h}^1$-atom supported in $B(o,r)$, } r\leq 1\}\lesssim e^{\rho \,t}.
\]
All throughout the proof, we let $\epsilon\coloneqq -b-d>0$, so that $m\in S^{-d-\epsilon}_\rho$. It will be crucial for the following to notice that by Lemma~\ref{Sobolev} (ii) with $\frac{1}{s}=\frac12-\left(-\frac{d+\epsilon}{n}\right)=1-\frac{1}{2n}+\frac{\epsilon}{n}$, H\"older inequality and the size condition of the atom we get
\begin{equation}\label{normaL2Tta}
\|\mathcal T_ta\|_{L^2} \lesssim \|a\|_{L^s} \lesssim\mu(B)^{-1+1/s} 	\lesssim r^{-\frac{1}{2}+\epsilon}.
\end{equation}

\bigskip

{\bf{{\emph{Case I: $t\geq 1/2$.}}}}  

Choose $J$ such that $J-2\leq  t+\frac{2}{10}\leq J-1$. Then for every $j\geq J$, the function $a\ast (\psi_jk_t)$ is supported in $B(o,j+1+r)$. By Lemma~\ref{aastgamma} and estimate~\eqref{Kttgrande} we obtain 
\[
\begin{aligned}
\|a\ast (\psi_jk_t)\|_{\mathfrak h^1}&\lesssim(\mu(B(o,j+r+1))^{1/2}\,\| \psi_jk_t \|_{L^2}\\
&\lesssim e^{\rho j}\,  \left(\int_{j-1}^{j+1}   e^{2\rho t}  e^{-4\rho s}\, |t-s|^{-4}\,e^{2\rho s}\,ds\right)^{1/2}\lesssim   e^{\rho t}   \,  |t-j|^{-2}.
\end{aligned}
\]
Thus
\begin{equation}\label{psijtgrande}
\begin{aligned}
\sum_{j=J}^{\infty} \|a\ast (\psi_jk_t)\|_{\mathfrak h^1}&\lesssim  e^{\rho t} \sum_{j=J}^{\infty}(j-t)^{-2}\lesssim  e^{\rho t} \int_J^{\infty}\, \frac{du}{(u-t)^2}\lesssim e^{\rho t},
\end{aligned}
\end{equation}
where we have used the fact that $J-2\leq  t+\frac{2}{10}\leq J-1$. 

\smallskip

{\bf{Subcase IA: $r\leq \frac{1}{10}$}}. 

Let $\phi_0$ be a smooth function taking values in $[0,1]$ supported in $B(o,3r)$ such that 
\[
\phi_0+\sum_{i=1}^{I_1}\phi_i+\sum_{i=I_1+1}^{I_2}\phi_i+\sum_{h=3}^{H_1}\eta_h+\sum_{h=3}^{H_2}\omega_h+\sum_{j=J}^{\infty}\psi_j=1 \qquad \mbox{in}\quad  \mathbb X\setminus A^{t+10r}_{t-10 r},
\]
where $\phi_i, \eta_h, \omega_h, \psi_j$ are defined by formulae~\eqref{phii},~\eqref{etahomegah},~\eqref{psij} and 
\[
\begin{aligned}
2^{I_1-1}r&\leq \tfrac{1}{10}\leq 2^{I_1+1}r,\\
2^{I_2-1}r &\leq  t-\tfrac{2}{10}\leq 2^{I_2+1}r, \\
t-2^{H_1+1} r &\leq t-\tfrac{2}{10} \leq t-2^{H_1-1} r,\\
t+2^{H_2-1} r&\leq t+\tfrac{2}{10} \leq t+2^{H_2+1}r.
\end{aligned}
\]
Define
\[
\sigma_t=\left[1-\phi_0+\sum_{i=1}^{I_1}\phi_i+\sum_{i=I_1+1}^{I_2}\phi_i+\sum_{h=3}^{H_1}\eta_h+\sum_{h=3}^{H_2}\omega_h+\sum_{j=J}^{\infty}\psi_j\right]k_t,
\]
so that
\[
\mathcal T_ta=a\ast (\phi_0k_t)+\sum_{i=1}^{I_2}a\ast (\phi_i k_t)+\sum_{h=3}^{H_1}a\ast (\eta_h k_t)+\sum_{h=3}^{H_2}a\ast (\omega_h k_t)+a\ast \sigma_t + \sum_{j=J}^{\infty}a\ast (\psi_j k_t).
\]
The $\mathfrak{h}^1$-norm of the last term of the sum has been already estimated in~\eqref{psijtgrande}. We now concentrate on the remaining terms.

The function $a\ast (\phi_0k_t)$ is supported in $B(o,4r)$ and by Lemma~\ref{normah1funzioneL2}
\begin{equation}\label{phi0kt}
\begin{aligned}
\|a\ast (\phi_0k_t)\|_{\mathfrak h^1}&\leq  \mu(B(o,4r))^{1/2} \|a\ast (\phi_0k_t)\|_{L^2} \lesssim r^{n/2} \| \mathcal T_t\|_{L^2\rightarrow L^2} \|a\|_{L^2} \lesssim 1,
\end{aligned}
\end{equation}
where we have used the size condition of the atom and the fact that the norm of the operator $f\mapsto f\ast (\phi_0k_t)$ on $L^2(\mathbb X)$ is bounded by the norm of $\mathcal T_t$ on $L^2(\mathbb X)$ (see e.g.~\cite[proof of Theorem 3.1]{MMV}).

Consider now the cases $i=1,\dots,  I_1$. The function $a\ast (\phi_ik_t)$ is supported in $B(o,(2^{i+1}+1)r)$. By Lemma~\ref{aastgamma} and by estimates~\eqref{Kttgrande} and~\eqref{Kt'tgrande} we obtain that 
\[
\begin{aligned}
\|a\ast (\phi_ik_t)\|_{\mathfrak h^1}&\lesssim (\mu(B(o,(2^{i+1}+1)r))^{1/2}\,r\,\|\nabla( \phi_ik_t) \|_{L^2}\\
&\lesssim(2^ir)^{n/2}\,r\,\left(\int_{ 2^{i-1} r}^{2^{i+1}r}   [  (2^ir)^{-2} s^{-2d-2+2\epsilon}+s^{-2d-4+2\epsilon}  ]     s^{n-1}\,ds\right)^{1/2}\\
&\lesssim (2^{i})^{\epsilon + (n-3)/2} r^{\epsilon + (n-1)/2}.
\end{aligned}
\]
Thus, since $I_1\asymp  \log_2\left(\frac{1}{10r}\right)$, we get
\begin{equation*}
\begin{aligned}
\sum_{i=1}^{I_1}\|a\ast (\phi_ik_t)\|_{\mathfrak h^1} &\lesssim  r^{\epsilon + (n-1)/2} \int_1^{\log_2\left(\frac{1}{10r}\right)}   (2^{u})^{\epsilon + (n-3)/2} \, du \lesssim r  .
\end{aligned}
\end{equation*}
Consider now the cases when $i=I_1+1,\dots,  I_2 $. The function $a\ast (\phi_ik_t)$ is supported in $B(o,(2^{i+1}+1)r)$.  By Lemma~\ref{aastgamma} and by estimates~\eqref{Kttgrande} and~\eqref{Kt'tgrande} we obtain that 
\[
\begin{aligned}
\|a\ast (\phi_i k_t)\|_{\mathfrak h^1}&\lesssim(\mu(B(o,(2^{i+1}+1)r))^{1/2}\,r\,\|\nabla( \phi_ik_t) \|_{L^2}\\
&\lesssim e^{\rho 2^ir} r\,\left(\int_{ 2^{i-1}r}^{2^ir}   [(2^ir)^{-2}  e^{-2\rho s}\, +e^{-2\rho s}]e^{2\rho s}\,ds\right)^{1/2} \lesssim e^{\rho 2^ir}\,    r^{3/2}\,2^{i/2}.
\end{aligned}
\]
Thus, since $2^{I_2} r\asymp t-\frac{1}{10}$, we get
\begin{equation*}
\begin{aligned}
\sum_{i=I_1+1}^{ I_2} \|a\ast (\phi_ik_t)\|_{\mathfrak h^1}&\lesssim r^{3/2}  \sum_{i=I_1+1}^{ I_2}   e^{\rho 2^ir}\,    2^{i/2}   \lesssim r^{3/2} \int_{  I_1+1}^{ I_2}  e^{\rho r2^u}  2^{u/2}\, du\\
&\lesssim     r^{3/2} \int_{2^{I_1+1}}^{2^{I_2}} v^{-1/2} e^{\rho v r}\, dv \lesssim r^{1/2}e^{\rho t}.
\end{aligned}
\end{equation*}
Consider now $3\leq h\leq H_1$. By the triangular inequality, the function $a\ast (\eta_h k_t)$ is supported in $A^{t-(2^{h-1}-1)r}_{t-(2^{h+1}+1)r}$, has vanishing integral and by Lemma~\ref{normah1funzioneL2medianulla}, by~\eqref{L2convgrad} and estimates~\eqref{Kttgrande} and~\eqref{Kt'tgrande}
\[
\begin{aligned}
\|a\ast (\eta_h k_t)\|_{\mathfrak h^1}&\lesssim \log(1/(2^h r))  (2^h r )^{1/2}\,e^{\rho t}\,r\,\|\nabla(\eta_h k_t)\|_{L^2}\\
&\lesssim e^{\rho t}r	\, (2^h r)^{-1+\epsilon }  \log(1/(2^h r))
\end{aligned}
\]
since $\|\nabla(\eta_h k_t)\|_{L^2}\lesssim (2^h r)^{-3/2+\epsilon} $. Using the fact that $2^{H_1}r\asymp \frac{2}{10}$ and then changing variables $2^h r=v$ we obtain
\begin{equation*}\label{etahkt}
\begin{aligned}
\sum_{h=3}^{H_1} \|a\ast (\eta_h k_t)\|_{\mathfrak h^1}
&\lesssim e^{\rho t} r \sum_{h=3}^{H_1} (2^h r)^{-1+\epsilon }  \log(1/(2^h r))\\
&\lesssim  e^{\rho t} r \int_{3}^{H_1} (2^h r)^{-1+\epsilon }  \log(1/(2^h r))\, dh \\
&\lesssim e^{\rho t} r  \int_{8r}^{2/10}  v^{-2+\epsilon } \log(1/v) \, dv \lesssim e^{\rho t} r^\epsilon \log(1/r) \lesssim e^{\rho t}. 
\end{aligned}
\end{equation*}
Similar computations can be done for $a* (\omega_h k_t)$, proving that
\[
\sum_{h=3}^{H_2} \| a* (\omega_h k_t)\|_{\mathfrak{h}^1}\lesssim e^{\rho t}.
\]
It remains to consider $a\ast \sigma_t$, where $\sigma_t$ is the singular part of the kernel supported in $A_{t-10r}^{t+10r}$. By the triangular inequality, the function $a\ast \sigma_t$ is supported in $A_{t-11r}^{t+11r}$. For every $x\in A_{t-11r}^{t+11r}$, we have 
\[
\mathcal T_ta(x)=a\ast\sigma_t(x)+a\ast (\eta_{3} k_t)(x)+a \ast  (\omega_{3} k_t)(x),
\]
so that
\begin{align}\label{normaL2aastsigmat}
\|a\ast\sigma_t\|_{L^2}&\leq \|\mathcal T_ta\|_{L^2}+\|a\ast (\eta_3 k_t)\|_{L^2}+\|a\ast (\omega_3 k_t)\|_{L^2}\nonumber\\
&\lesssim r^{-1/2+\epsilon }+r\|\nabla(\eta_3k_t)\|_{L^2}+r\|\nabla(\eta_3 k_t)\|_{L^2} \lesssim r^{-1/2+\epsilon}.
\end{align}
The second inequality follows from~\eqref{normaL2Tta} and~\eqref{L2convgrad}, while the third follows from the computations we made before for $\nabla (\eta_h k_t)$ and a similar computation for $\nabla  (\omega_h k_t)$.

We deduce from Lemma~\ref{normah1funzioneL2medianulla} and~\eqref{normaL2aastsigmat} that
\begin{equation*}
\|a\ast\sigma_t\|_{\mathfrak h^1} \lesssim \log(1/r) e^{\rho t} r^{1/2}r^{-1/2+\epsilon}\lesssim e^{\rho t}.
\end{equation*}

{\bf{Subcase IB: $\frac{1}{10}<r\leq 1$.}}
 
Choose two smooth cutoff functions $\phi_0$ and $\phi_t$ with values in $[0,1]$ such that 
\[
\begin{aligned}
&\mathrm{supp}\,(\phi_0)\subseteq B(o,3),\qquad \mathrm{supp}\,(\phi_t)\subseteq A^{t-\frac{1}{10}}_2\\
&\phi_0+\phi_t+\sum_{j=J}^{\infty}\psi_j=1\qquad \mbox{in} \quad \mathbb X\setminus A^{t+\frac{1}{10}}_{t-\frac{1}{10}},
\end{aligned}
\]
(if $t-1/10<2$, then just $\phi_t\equiv 0$) and define 
\[
\sigma_t=\left[1-\phi_0-\phi_t-\sum_{j=J}^{\infty}\psi_j\right]k_t.
\]
The convolution of $a$ with the sum of the $\psi_j$'s has been already estimated in~\eqref{psijtgrande}. The function $a\ast (\phi_0k_t)$ is supported in $B(o,3+r)$ and by Lemma~\ref{normah1funzioneL2} 
\begin{equation}\label{phi0kt2}
\|a\ast (\phi_0k_t)\|_{\mathfrak h^1}\lesssim  \mu(B(o,4))^{1/2}\,\|\mathcal T_t\|_{L^2\rightarrow L^2} \,\|a\|_{L^2} \lesssim 1 ,
\end{equation}
where we argued as in~\eqref{phi0kt}. By Lemma~\ref{aastgamma} and estimates~\eqref{Kttgrande} we get 
\begin{equation}\label{phitkt}
\begin{aligned}
\|a\ast (\phi_tk_t)\|_{\mathfrak h^1}&\lesssim \mu (B (o,t-\tfrac{1}{5}+r ) )^{1/2}\,\|\phi_tk_t\|_{L^2}\\
&\lesssim e^{\rho t}\,\left(  \int_2^{t-\frac{1}{5}} e^{-2\rho s}|t-s|^{-4} e^{2\rho s}\, ds \right)^{1/2}\lesssim e^{\rho t}.
\end{aligned}
\end{equation}
It remains to estimate the $\mathfrak h^1$-norm of $a\ast \sigma_t$, which is supported in $A^{t+r +1/10}_{t-r -1/10}$. Since
\[
\mathcal T_ta(x)=a\ast\sigma_t(x)+a\ast (\phi_tk_t)(x) +a\ast (\psi_Jk_t)(x) \qquad \forall \; x\in A^{t+r+ 1/10}_{t-r- 1/10},
\]
then
\[
\begin{aligned}
\|a\ast\sigma_t\|_{L^2}&\leq \|\mathcal T_ta\|_{L^2}+\|a\ast (\phi_tk_t)\|_{L^2}+\|a\ast (\psi_Jk_t)\|_{L^2}\\
&\leq \|\mathcal T_t\|_{L^2\rightarrow L^2}\,\|a\|_{L^2}\,+\|\phi_tk_t\|_{L^2}+\| \psi_Jk_t\|_{L^2}\lesssim 1,
\end{aligned}
\]
which follows from~\eqref{normaL2Tta} and the computations we made in~\eqref{phitkt} and~\eqref{psijtgrande}. Thus
\begin{equation*}\label{sigmat2}
\|a\ast\sigma_t\|_{\mathfrak h^1}\lesssim \mu\left( B(o,t+\tfrac{1}{10}+r )\right)^{1/2}\,\|a\ast\sigma_t\|_{L^2}\lesssim e^{\rho t}.
\end{equation*}
The proof in the case $t\geq 1/2$ is then complete.

\bigskip

{\bf{{\emph{Case II: $t<1/2$.}}  }} 

For every $j\geq 2$ by Lemma~\ref{aastgamma} and estimates~\eqref{Kttpiccolo} we get
\[\begin{aligned}
\|a\ast(\psi_jk_t)\|_{\mathfrak h^1}&\lesssim \mu(B(o,j+2))^{1/2}\,\|\psi_jk_t\|_{L^2}\\
&\lesssim  e^{\rho j}\left( \int_{j-1}^{j+1}e^{-4\rho s}(1+|t-s|^{-2})^{2} e^{2\rho s}	\, ds\right)^{1/2}\lesssim (j-t)^{-2}\lesssim j^{-2},
\end{aligned}
\]
where the functions $\psi_j$ are defined in~\eqref{psij}. Thus
\begin{equation}\label{psijtpiccolo}
\sum_{j=2}^{\infty}\|a\ast(\psi_jk_t)\|_{\mathfrak h^1}\lesssim \sum_{j=2}^{\infty}j^{-2}\lesssim 1.
\end{equation}

{\bf{{Subcase IIA: $r\leq \frac{t}{20}$. }  }}

Let $\phi_0$ be a cutoff function supported in $B(o,3r)$ taking values in $[0,1]$ such that
\[
\begin{aligned}
&\phi_0+\sum_{i=2}^I\phi_i+\sum_{i=I_1}^{I_2}\phi_i+   \sum_{j=2}^{\infty}\psi_j=1\qquad \mbox{in}\quad \mathbb X\setminus A^{t+10 r}_{t-10r},
\end{aligned}
\]
where the $\phi_i$'s are defined by~\eqref{phii} and 
\[
\begin{aligned}
2^{I-1}r&<t-10r<2^{I+1}r,\\
2^{I_1-1}r&<t+10 r<2^{I_1+1}r,\\
2^{I_2-1}r&<1<2^{I_2+1}r.
\end{aligned}
\] Define 
\[
\sigma_t=\left[1-\phi_0-\sum_{i=3}^I \phi_i-\sum_{i=I_1}^{I_2}\phi_i-\sum_{j=1}^{\infty}\psi_j\right]\,k_t.
\]
The $\mathfrak{h}^1$-norm of the convolution with the $\psi_j$'s has been already estimated in~\eqref{psijtpiccolo}. Since $a\ast(\phi_0k_t)$ is supported in $B(o,4r)$
\begin{equation*}
\|a\ast(\phi_0 k_t)\|_{\mathfrak h^1}\lesssim \mu((B(o,4r)))^{1/2}  \|a\ast (\phi_0 k_t)\|_{L^2}\lesssim r^{n/2}\|a\|_{L^2}\,\|\mathcal T_t\|_{L^2\rightarrow L^2}\lesssim  1 ,
\end{equation*}
where we argued as in~\eqref{phi0kt}. For every $i\in\{2,\dots,I\}\cup\{I_1,\dots,I_2\}$, the function $a\ast (\phi_ik_t)$ is supported in $B(o,2^{i+1}r+r)$ and by Lemma~\ref{aastgamma} and estimates~\eqref{Kttpiccolo} and~\eqref{Kt'tpiccolo}
\[
\begin{aligned}
\|a\ast (\phi_ik_t)\|_{\mathfrak h^1}&\lesssim  \, \mu(B(o,2^{i+1}r+r))^{1/2}\,r\,\|\nabla (\phi_i k_t)\|_{L^2}\\
&\lesssim (2^ir)^{n/2}\,r\,\Bigg( \int_{2^{i-1}r}^{2^{i+1}r} [ (2^ir)^{-2}s^{-2-2d+2\epsilon}+ (2^ir)^{-2}s^{-2d}|t-s|^{-2+2\epsilon}\\&\quad+ s^{-2d-4+2\epsilon}+ s^{-2d}|t-s|^{-4+2\epsilon} +s^{-2d-2}|t-s|^{-2+2\epsilon}    ]s^{n-1}\,ds\Bigg)^{1/2}\\
&\lesssim r\, \left[  (2^ir)^{\frac{n-3}{2}+\epsilon}+(2^ir)^{\frac{n-1}{2}}|t-2^ir|^{-1+\epsilon}+(2^ir)^{\frac{n+1}{2}} |t-2^ir|^{-2+\epsilon}    \right].
\end{aligned}
\]
Thus
\[
\begin{aligned}
\sum_{i=2}^I\|a\ast (\phi_ik_t)\|_{\mathfrak h^1}&\lesssim r^{\frac{n-1}{2}+\epsilon }\, \int_{2}^I(2^u)^{\frac{n-3}{2}+\epsilon}\, du+r^{\frac{n+1}{2}}\, \int_{2}^I(2^u)^{\frac{n-1}{2}}|t-2^u r|^{-1+\epsilon}\, du\\ &\quad +r^{\frac{n+3}{2}}\, \int_{2}^I(2^u)^{\frac{n+1}{2}}|t-2^u r|^{-2+\epsilon}\,du.
\end{aligned}
\]
By the change of variables $2^u r=w$ and recalling that $2^Ir \asymp t-10 r<t<1/2$,
\begin{align*}\label{phiitpiccolo}
\sum_{i=3}^I\|a\ast (\phi_ik_t)\|_{\mathfrak h^1}& \lesssim r +r \int_{4r}^{2^Ir} w^{\frac{n-3}{2}}|t-w|^{-1+\epsilon}\, dw + r \int_{4r}^{2^Ir} w^{\frac{n-1}{2}}|t-w|^{-2+\epsilon}\,dw \\ & \lesssim r^\epsilon \lesssim 1
\end{align*}
since $|t-w|\geq |t-4r|\gtrsim r$. Arguing as before, we can also prove that  
\begin{equation*}\label{phiitpiccolo2}
\sum_{i=I_1}^{I_2}\|a\ast (\phi_ik_t)\|_{\mathfrak h^1}\lesssim 1.  
\end{equation*}
It remains to consider $a\ast \sigma_t$, where $\sigma_t$ is the singular part of the kernel supported in $A_{t-10r}^{t+10r}$. By the triangular inequality, $a\ast \sigma_t$ is supported in $A_{t-11r}^{t+11r}$. For every $x\in A_{t-11r}^{t+11r}$, we have 
\[
\mathcal T_ta(x)=a\ast\sigma_t(x)+a\ast (\phi_I k_t)(x)+a \ast  (\phi_{I_1} k_t)(x),
\]
so that
\[
\begin{aligned}
\|a\ast\sigma_t\|_{L^2}&\leq \|\mathcal T_ta\|_{L^2}+\|a\ast (\phi_Ik_t)\|_{L^2}+\|a\ast (\phi_{I_1} k_t)\|_{L^2}\\
&\lesssim r^{-1/2+\varepsilon}+r\|\nabla(\phi_Ik_t)\|_{L^2}+r\|\nabla(\phi_{I_1} k_t)\|_{L^2}\lesssim r^{-1/2+\epsilon},
\end{aligned}
\]
where we have applied~\eqref{normaL2Tta} and the computations we made above. Then by Lemma~\ref{normah1funzioneL2medianulla}
\begin{equation*}\label{sigmattpiccolo}
\|a\ast\sigma_t\|_{\mathfrak h^1}\lesssim  \log(1/r)\mu(A_{t-11r}^{t+11r})^{1/2} \|a\ast\sigma_t\|_{L^2}\lesssim \log(1/r) r^{\epsilon}\lesssim 1 .
\end{equation*}

{\bf{{Subcase IIB: $ \frac{t}{20}<r\leq 1$.}  }}
 
Notice that $t+10r<30r$. We choose a smooth cutoff function $\phi_0$ supported  in $B(o,30r)$ taking values in $[0,1]$ such that
\[
\begin{aligned}
&\phi_0+\sum_{i=5}^I\phi_i+\sum_{j=2}^{\infty}\psi_j=1\
\end{aligned}
\]
in $\mathbb{X}$, where $I$ is such that $2^{I-1}r<1<2^{I+1}r$. We split the kernel $k_t$ accordingly as we did before. Then $a\ast (\phi_0k_t)$ is supported in $B(o,31r)$ and
\begin{equation*}
\|a\ast(\phi_0 k_t)\|_{\mathfrak h^1}\lesssim \mu(B(o,31r))^{1/2} \|a\|_{L^2}\,\| \mathcal T_t\|_{L^2\rightarrow L^2}\lesssim 1,
\end{equation*}
where we argued as in~\eqref{phi0kt}. For every $i=5,\dots,I$ by Lemma~\ref{aastgamma} and estimate~\eqref{Kttpiccolo} one can see that
\begin{equation*}
\begin{aligned}
\|a\ast (\phi_ik_t)\|_{\mathfrak h^1}&\lesssim (2^i r)^{n/2} \| \phi_ik_t\|_{L^2}
\lesssim (2^ir)^{\frac{n-1}{2}+\epsilon} + (2^i r)^{\frac{n+1}{2}} |t-2^i r|^{-1+\epsilon},
\end{aligned}
\end{equation*}
which yields 
\begin{align*}
\sum_{i=5}^I \|a\ast (\phi_ik_t)\|_{\mathfrak h^1} \lesssim r^{\frac{n-1}{2}+\epsilon} \int_{2^5}^{2^I} v^{\frac{n-3}{2}+\epsilon}\, dv  + \int_{32r}^1 v^{\frac{n-1}{2}}|t-v|^{-1+\epsilon}\, dv  \lesssim 1
\end{align*}
where we used the fact that $2^{I-1}r<1<2^{I+1}r$. This concludes the proof of the case $t<1/2$ and of the theorem. \hfill\ensuremath{\blacksquare}

\bigskip

{\bf{Acknowledgments.}} 	The authors are members of the Gruppo Nazionale per l'Analisi Matematica, la Probabilit\`a e le loro Applicazioni (GNAMPA) of the Istituto Nazionale di Alta Matematica (INdAM).
This work was partially supported by the Progetto PRIN 2015 ``Variet\`a reali e complesse: geometria, topologia e analisi armonica''.

The authors would like to thank Stefano Meda, Fulvio Ricci and Peter Sj\"ogren for helpful discussions about this work.

\end{document}